\documentclass{ifacconf}

\usepackage{graphicx}      
\usepackage{natbib}        
\usepackage{graphics} 
\usepackage[dvpinames]{xcolor}
\usepackage{microtype}
\usepackage{amsmath} 
\usepackage{amssymb}  

\usepackage{algorithm}
\usepackage{etoolbox}
\usepackage{todonotes}

\let\classAND\AND
\let\AND\relax
\usepackage{algorithmic}

\let\AND\classAND
\AtBeginEnvironment{algorithmic}{\let\AND\algoAND}

\floatname{algorithm}{Algorithm}
\renewcommand{\algorithmicrequire}{\textbf{Input:}}
\renewcommand{\algorithmicensure}{\textbf{Output:}}

\usepackage{xspace}

\def\IR{{\mathbb R}}
\def\IC{{\mathbb C}}

\def\IL{{\mathbb L}}
\def\IM{{\mathbb M}}

\def\IV{{\mathbb V}}
\def\IW{{\mathbb W}}

\newcommand{\bA}{{\bf A}}
\newcommand{\bB}{{\bf B}}
\newcommand{\bC}{{\bf C}}

\newcommand{\bE}{{\bf E}}

\newcommand{\bL}{{\bf L}}

\newcommand{\bH}{{\bf H}}
\newcommand{\bP}{{\bf P}}

\newcommand{\bQ}{{\bf Q}}
\newcommand{\bW}{{\bf W}}
\newcommand{\bR}{{\bf R}}

\newcommand{\bX}{{\bf X}}
\newcommand{\bx}{{\bf x}}
\newcommand{\by}{{\bf y}}
\newcommand{\bu}{{\bf u}}
\newcommand{\bV}{{\bf V}}
\newcommand{\bU}{{\bf U}}

\newcommand{\bfe}{{\bf e}}

\newcommand{\bxi}{ \boldsymbol{\xi} }
\newcommand{\bchi}{ \boldsymbol{\chi} }

\newcommand{\bfz}{{\mathbf 0}}

\newcommand{\cK}{ {\cal K} }

\newcommand{\bPhi}{ \boldsymbol{\Phi} }

\newcommand{\Htwo}{{{\cal H}_2}}
\newcommand{\Hinf}{{{\cal H}_{\infty}}}

\newcommand{\nq}{{{}N_q}}

\newcommand{\quadW}{\widetilde{\bW}}

\newcommand{\quadQ}{\widetilde{\bQ}}

\newcommand{\quadL}{\widetilde{\bL}}

\newcommand{\quadLL}{\widetilde{\IL}}
\newcommand{\quadMM}{\widetilde{\IM}}
\newcommand{\quadVV}{\widetilde{\IV}}
\newcommand{\quadWW}{\widetilde{\IW}}

\newcommand{\LF}{\textsf{LF}\xspace}
\newcommand{\RK}{\textsf{IRKA}\xspace}
\newcommand{\SK}{\textsf{ISRK}\xspace}
\newcommand{\QSK}{\textsf{Quad-ISRK}\xspace}

\newcommand{\QBT}{\textsf{QuadBT}\xspace}
\newcommand{\BT}{\textsf{BT}\xspace}
\newcommand{\imunit}{{\dot{\imath\hspace*{-0.2em}\imath}}}

\newcommand{\bAtil}{\widetilde{\bf A}}
\newcommand{\bEtil}{\widetilde{\bf E}}
\newcommand{\bBtil}{\widetilde{\bf B}}
\newcommand{\bCtil}{\widetilde{\bf C}}

\newcommand{\bomega}{{\bf \omega}}

\usepackage{bbm}
\renewcommand{\algorithmicrequire}{\textbf{Input:}}
\renewcommand{\algorithmicensure}{\textbf{Output:}}

\newtheorem{theorem}{Theorem}[section]

\newtheorem{lemma}[theorem]{Lemma}
\newtheorem{corollary}[theorem]{Corollary}
\newtheorem{remark}[theorem]{Remark}

\newtheorem{definition}[theorem]{Definition}

\begin{document}
\begin{frontmatter}

\title{A non-intrusive data-based reformulation \\ of a  hybrid projection-based \\ model reduction method\thanksref{footnoteinfo}} 

\thanks[footnoteinfo]{
 The work of S.~Gugercin
was supported in parts by the National Science Foundation under Grant No. DMS-2318880.}

\author[First]{Ion Victor Gosea} 
\author[Second]{Serkan Gugercin} 
\author[Third]{Christopher Beattie}

\address[First]{
        Max Planck Institute for Dynamics of Complex Technical Systems, Magdeburg, Germany (e-mail: gosea@mpi-magdeburg.mpg.de).}
\address[Second]{Department of Mathematics, Virginia Tech, Blacksburg, VA, USA (e-mail: gugercin@vt.edu)}
\address[Third]{Department of Mathematics, Virginia Tech, Blacksburg, VA, USA (e-mail: beattie@vt.edu)}

\begin{abstract}                
    We present a novel data-driven reformulation of the iterative SVD-rational Krylov algorithm (\SK), in its original formulation
    a Petrov-Galerkin (two-sided) projection-based iterative method for model reduction combining rational Krylov subspaces (on one side) with Gramian/SVD based subspaces (on the other side). We show that in each step of \SK, we do not necessarily require access to the original system matrices, but only to input/output data in the form of the system's transfer function, evaluated at particular values (frequencies). Numerical examples illustrate the efficiency of the new data-driven formulation. 
\end{abstract}

\begin{keyword}
Efficient strategies for large-scale complex systems, Linear systems, Identification and model reduction, Time-invariant systems, Data-based methods, Transfer function.
\end{keyword}

\end{frontmatter}

\section{Introduction}

Model order reduction (MOR) techniques aim to approximate large-scale, complex models of dynamical systems, with much simpler and smaller models that can be easily simulated, controlled, or further analyzed. A generic goal for most MOR methods is to compute reduced-order models (ROMs) that are close, in terms of response, characteristics, and properties, to the original complex,  large-scale model. 
Additionally, with the ever-increasing need to incorporate measured data in the complexity reduction, identification, or control design procedures, it is of interest to develop data-driven methods. These methods represent an alternative approach to classical (intrusive) methods, which typically rely on explicit access to a large-scale model, e.g., an explicit state-space representation. Unlike the intrusive methods, data-driven reduced-order modeling does not require explicit knowledge of the model's structure or matrices. Instead, low-order models are constructed directly from measured data (either in the time or frequency domains).

System-theoretical MOR methods are the focus of the current work, with balanced truncation (\BT) in \cite{moore1981principal,mullis1976synthesis}
and interpolatory approaches~\cite{AntBG20} being relevant methodologies. \BT remains a standard for system-theoretical MOR methodologies, decades after its original development. Its main attributes are a priori error bounds, stability preservation, and robust numerical implementations (that could cope with the large-scale setup and with several extensions of the method to different structures or constraints). We refer the reader to  \cite{benner2017model} which covers recent extensions and provides relevant insights into the method's development. In the second class of methods, based on interpolation theory, the original system function transfer function and its derivatives are matched at selected points. This yields advantages through the use of (rational) Krylov subspace methods, which are easy to implement using readily available algorithms from the field of numerical linear algebra~\cite{AntBG20}. The iterative rational Krylov algorithm (\RK) \cite{gugercin2008h_2}  and the realization-free methods such as the Loewner framework (\LF) in \cite{mayo2007fsg} are two examples that are relevant to the current work. \RK provides an iterative approach for finding optimal ROMs in the so-called $\mathcal{H}_2$ norm and shows how to find optimal interpolation points in the $\mathcal{H}_2$ norm. \LF computes ROMs by arranging measured data in Loewner matrices and extracting dominant features.

Recently, a data-driven quadrature-based balanced truncation method (\QBT) was proposed in \cite{gosea2022data}. This approach makes use of the quadrature representation of infinite Gramians only implicitly and circumvents the need to access the system matrices explicitly, making \BT non-intrusive. Such data-driven formulations are of particular importance whenever a state-space form of a model is not available as they utilize frequency-domain measurements, i.e., samples of the underlying system's transfer function, to compute ROMs. Extensions of \QBT were recently proposed in \cite{morReiGG24}.

In this paper, we propose a non-intrusive data-based reformulation of an iterative, projection-based (double-sided) model reduction method originally proposed in \cite{gugercin2008iterative}. The original method under consideration, the iterative SVD-rational Krylov (\SK) is based on a Petrov-Galerkin projection and combines rational Krylov-based and balancing (Gramians) techniques. The latter part depends on the observability  Gramian, $\mathbf{Q}$, while the rational Krylov side is obtained via iteratively corrected interpolation points. It was shown in the original contribution \cite{gugercin2008iterative}, that the reduced-order model thus obtained is asymptotically stable, matches interpolation at selected points, and solves a restricted $\mathcal{H}_2$ minimization problem.
The reformulation of \SK proposed here, referred to as \QSK, is inspired by the \QBT method~\cite{gosea2022data}, in the sense that it makes (implicit) use of the quadrature-based representation of the $\mathbf{Q}$-Gramian. Then, at each step of the iteration, reduced-order matrices are put together in terms of transfer function measurements; the sampling points are computed as the mirrored poles of the model constructed at the previous iteration.

The paper is organized as follows: after the introduction has been set up,  Section \ref{sec:prel} gathers various definitions, and briefly introduces the system Gramians, together with the iterative projection approach in \cite{gugercin2008iterative}. Then, Section \ref{sec:main} provides the main contribution of the paper, i.e., a data-driven reformulation of the latter method, which can be computed without having access to the system's realization matrices, but using only transfer function evaluations. Then, in Section \ref{sec:num} we illustrate the applicability of the method to a numerical test case: a classical MOR benchmark example. Finally, Section \ref{sec:conc} presents the conclusions, and provides an outlook to possible future endeavors. 


\section{Preliminaries}
\label{sec:prel}

In this work, we consider linear time-invariant (LTI) systems described in state-space by the  equations
\begin{align}
\begin{split} \label{OrigSys}
\bE\, \dot{\bx}(t) &= \bA\, \bx(t) + \bB\, \bu(t), \\
\by(t) &= \bC \bx(t),
\end{split}
\end{align}
where the input mapping is given by $\bu : \mathbb{R} \rightarrow \mathbb{R}^m$, the (generalized) state trajectory/variable is
$\bx: \mathbb{R} \rightarrow \mathbb{R}^n$, and the output mapping is $y : \mathbb{R} \rightarrow \mathbb{R}^p$. The system matrices are given by $\bE,\,\bA \in \mathbb{R}^{n \times n}$, and  $\bB\in \mathbb{R}^{n\times m}$, $\bC \in \mathbb{R}^{p\times n}$. We assume that the matrix $\bE$ is nonsingular and that the matrix pencil $(\bA,\bE)$ is \emph{asymptotically stable}, i.e., the eigenvalues of the matrix pencil $\lambda \bE - \bA$ are located in the complex open left half-plane.
The transfer function of the system \eqref{OrigSys} is a $p\times m$ rational matrix function, defined as
\begin{equation}\label{TF_def}
\bH(s)=\bC (s \bE -\bA)^{-1} \bB.    
\end{equation}
MOR aims to approximate the LTI system in (\ref{OrigSys}) with an LTI system of lower dimension, denoted with $r$, described in state-space by the  equations
\begin{align}
\begin{split} \label{RedSys}
\bE_r\, \dot{\bx}_r(t) &= \bA_r\, \bx_r(t) + \bB_r\, \bu(t), \\
\by_r(t) &= \bC_r \bx_r(t),
\end{split}
\end{align}
where the state variable is
$\bx_r: \mathbb{R} \rightarrow \mathbb{R}^r$;  the output mapping is $\by_r : \mathbb{R} \rightarrow \mathbb{R}^p$; and $\bE_r,\,\bA_r \in \mathbb{R}^{r \times r}$, and  $\bB\in \mathbb{R}^{r\times m}$, $\bC \in \mathbb{R}^{p\times r}$. The goal of MOR is that  $\by_r(t)$  approximates $\by(t)$ well for admissible inputs $\bu: \mathbb{R} \rightarrow \mathbb{R}^m$ (with bounded energy, i.e., $\Vert \bu \Vert_2 < \infty$). The transfer function of ROM~\eqref{RedSys} is the degree-$r$ rational function 
\begin{equation}\label{TFred_def}
\bH_r(s)=\bC_r (s \bE_r -\bA_r)^{-1} \bB_r.
\end{equation}
One can then use, e.g.,  the $\mathcal{H}_\infty$ or $\mathcal{H}_2$, distances between $\bH(s)$ and $\bH_r(s)$ to check the accuracy of the ROM; see, e.g., \cite{ACA05,AntBG20}.

When the state-space formulation as in~\eqref{OrigSys} is available, the intrusive Petrov-Galerkin projection is the common approach to form the reduced model~\eqref{RedSys}: In this case, one computes two model reduction bases $\bV_r,\bW_r \in \mathbb{C}^{n\times r}$ such that the ROM quantities in~\eqref{RedSys} are given by
\begin{align} \label{eq:proj}
\begin{split} 
\bE_r = \bW_r^T \bE \bV_r,
\, \quad & \bA_r = \bW_r^T \bA \bV_r, \\
\bB_r = \bW_r^T \bB, ~\mbox{and}~ &\bC_r = \bW_r^T \bC.
\end{split}
\end{align}

To simplify the notation, we focus on single-input single-output (SISO) LTI systems. Hence, for the rest of the paper, we assume $m = p =1$, i.e., $\bB$ and $\bC^T$ are column vectors.

\subsection{The system Gramians}
\label{sec:BT}

Two important quantities in systems theory are the reachability and observability Gramians. The reachability Gramian (typically denoted by $\bP$) provides a measure of how easily a state can be accessed from the zero state. The matrix $\bP$ is expressed in the frequency domain as:
\begin{align} \label{gram_P_freq}
\begin{split}
\hspace{-2mm}\bP =
\frac{1}{2\pi}\int_{-\infty}^{\infty} (\imunit \zeta\bE -\bA)^{-1} \bB \bB^T (-\imunit \zeta\bE^T -\bA^T)^{-1} d \zeta,
\end{split}
\end{align}
where $\imunit$, as $\imunit^2 = -1$, represents the imaginary unit. An alternative definition is possible in the time domain (see \cite{ACA05}).

Additionally, the observability Gramian describes how easily a state
can be observed, and is represented as $\bE^T\bQ\bE$ where $\bQ$ is defined in the frequency domain as 
\begin{align} \label{gram_Q}
\begin{split}
\hspace{-3mm} \bQ = \frac{1}{2\pi}\int_{-\infty}^{\infty} (-\imunit \omega\bE^T -\bA^T)^{-1} \bC^T \bC (\imunit \omega\bE -\bA)^{-1} d \omega.
\end{split}
\end{align}
 The Gramians $\bP$ and $\bQ$ satisfy the following generalized Lyapunov equations:
\begin{align}  \label{lyapforP}
\bA \bP\bE^T + \bE \bP \bA^T + \bB\bB^T = \bfz, \\
\bA^T \bQ \bE + \bE^T \bQ \bA + \bC^T \bC = \bfz. \label{lyapforQ}
\end{align}
We emphasize that $\bE^T\bQ\bE$ is the observability Gramian in this scenario.  
The solutions, $\bP$ and $\bQ$, to \eqref{lyapforP} and \eqref{lyapforQ}, respectively, are symmetric positive-semidefinite matrices; if the original system as specified in \eqref{OrigSys} is \emph{minimal} then both matrices will be positive-definite. The so-called square-roots factors $\bL, \bU \in \mathbb{R}^{n \times n}$ can be computed (e.g., by means of a classical Cholesky factorization) so that
\begin{equation} \label{lyapFact}
\bP = \bU \bU^T \quad \mbox{and} \quad \ \bQ = \bL \bL^T.
\end{equation}
The  matrices $\bL$ and $\bU$, are essential in the ``square-root" implementation of \BT; see, e.g.,~\cite{ACA05}.

\subsection{The Iterative SVD-Rational Krlyov Method (\SK)}
\label{sec:ISRK}

\BT~construct the Petrov-Galerkin model reduction bases
$\bW_r$ and $\bV_r$ in~\eqref{eq:proj} using the Gramians. On the other hand, interpolatory methods (as we outline below) construct the bases to enforce rational interpolation of $\bH(s)$. Indeed, interpolation (at selected points) forms the necessary conditions for MOR in $\mathcal{H}_2$ norm~\cite{gugercin2008h_2}. Then, 
inspired by these $\Htwo$-optimal interpolatory conditions and \RK and Gramian-based model reduction,
the main idea behind \SK~\cite{gugercin2008iterative} is to combine 
these two frameworks where the model reduction basis $\bV_r$
is interpolation based and the basis $\bW_r$ is based on the observability Gramian. Before we go into a more detailed exposition of the method, we need to first set up the stage with notations and relevant concepts.

\begin{definition}
For matrices $\bA,\bB,\bE$ as in (\ref{OrigSys}), a positive integer $r$, and a vector of complex shift scalars $\bchi = \{\chi_1,\chi_2,\ldots,\chi_r\}$, which do \emph{not} overlap with the eigenvalues of the matrix pencil $\lambda \bE - \bA$,
we define \emph{the  primitive rational Krylov matrix}
$\cK(\bA,\bE,\bB,\bchi) \in \IC^{n \times r}$ as
\begin{align} \label{eq:K}
\begin{split} 
\cK(&\bA,\bE,\bB,\bchi) \\ &= \begin{bmatrix}
        (\chi_1 \bE -\bA)^{-1} \bB &  \cdots &  (\chi_r \bE -\bA)^{-1} \bB
    \end{bmatrix}.
    \end{split}
\end{align}
\end{definition}
Note that that the matrix $\cK(\bA,\bE,\bB,\bchi)$ is the unique solution of the generalized Sylvester equation  
\begin{equation}
    \bA \cK + \bB \bR = \bE \cK \bX,
\end{equation}
where $\bX = \text{diag}(\chi_1,\chi_2,\ldots,\chi_r) \in \IC^{r \times r}$, $\bR = \begin{bmatrix} 1 & 1 & \cdots & 1 \end{bmatrix}$ $\in \IC^{1 \times r}$, and the notation $\cK$ was used instead of $ \cK(\bA,\bE,\bB,\bchi)$. 
In practice, 
the matrix $\cK(\bA,\bE,\bB,\bxi)$ is seldom computed explicitly as in~\eqref{eq:K}. Instead, an orthonormal basis spanned by the columns of $\cK(\bA,\bE,\bB,\bxi)$ is computed.

In interpolatory MOR methods, typically the basis $\bV_r \in \IC^{n \times r}$ in~\eqref{eq:proj}  is chosen as
\begin{equation} \label{eq:Vr}
   \text{Ran}(\bV_r) =  \text{Orth}(\cK(\bA,\bE,\bB,\bchi)),
\end{equation}
where  \text{Orth} stands for an orthonormal basis. 
Additionally, in ~\eqref{eq:proj}, one typically constructs  $\bW_r \in \IC^{n \times r}$ as
\begin{equation} \label{eq:Wr}
   \text{Ran}(\bW_r) =  \text{Span}(\cK(\bA^T,\bE^T,\bC^T,\bxi)),
\end{equation}
where $\bxi = \{\xi_1,\xi_2,\ldots,\xi_r\}$ is the vector of complex shift scalars. Then, $\bH_r(s)$, the transfer of ROM constructed as in\eqref{eq:proj} interpolates $\bH(s)$ at the selected points, i.e.,
$$
\bH(\chi_i) = \bH_r(\chi_i)~~\mbox{and}~~
\bH(\xi_i) = \bH_r(\xi_i)~\mbox{for}~ i=1,2,\ldots,r.
$$
Higher-order interpolation can be achieved similarly. Indeed, optimal approximation in the $\mathcal{H}_2$ norm requires interpolating both $\bH(s)$ and $\bH'(s)$ at the mirror images of the poles of $\bH_r(s)$; see~\cite{AntBG20,gugercin2008h_2,MeiL67}. Moreover, the interpolation property above assumes direct solves for the shifted linear systems in forming~\eqref{eq:Vr} and \eqref{eq:Wr}. For the analysis of iterative methods in the interpolatory MOR setting, see~\cite{BGW12}.

In \cite{gugercin2008iterative}, an alternative way to choosing~$\bW_r$ was proposed, namely it was proposed to choose  $\bW_r$ as $\bW_r = \bQ \bE \bV_r$ where 
$\bV_r$ is given as in~\eqref{eq:Vr} and  $\bQ$ is the observability Gramian  as in~\eqref{lyapforQ}. 
Thus, while $\bV_r$ carries the interpolation property, $\bW_r$ also brings in information from the observability Gramian; therefore combining interpolatory projections with Gramian-based methods. 
Moreover, inspired by \RK and the interpolatory $\mathcal{H}_2$ optimality conditions, the \SK procedure is iterative, in the sense that the shifts in forming $\bV_r$, are updated at every step 
as the mirror images of the current reduced-order poles. The algorithm terminates once the shifts (interpolation points) converge. Basically, $\bV_r^{(k)}$ at step $k$ of the iteration in \SK, is given by
\begin{align}\label{Vrk_alg1}
  \hspace{-2mm}  \bV_r^{(k)} &= \begin{bmatrix}
        (\eta_1^{(k)} \bE -\bA)^{-1} \bB &  \cdots &  (\eta_r^{(k)} \bE -\bA)^{-1} \bB 
        \end{bmatrix},
\end{align}
with shifts (interpolation points) $\{\eta_1^{(k)},\eta_2^{(k)},\ldots,\eta_r^{(k)}\}$ at step $k$, and the $\bW_r^{(k)}$ at step $k$ of the iteration in \SK is 
\begin{align}\label{Wrk_alg1}
        \bW_r^{(k)} &=  \bQ \bE \bV_r^{(k)}.
\end{align}
Then, the ROM matrices at step $k$ are computed as:
\begin{align} \label{red_mat_stepk_alg1}
\begin{split}
    \bA_r^{(k)} &= \left[\bW_r^{(k)}\right]^T \bA \bV_r^{(k)}, \ \ \bE_r^{(k)} = \left[\bW_r^{(k)}\right]^T \bE \bV_r^{(k)},\\
        \bB_r^{(k)} &= \left[\bW_r^{(k)}\right]^T \bB, \ \ \bC_r^{(k)} =  \bC \bV_r^{(k)}.
        \end{split}
\end{align}
Then, in the next step, the shifts (interpolation points) are updated as 
$$
\eta_i^{(k+1)} \leftarrow -\tt{eig}(\bA_r^{(k)},   \bE_r^{(k)}),
$$
where $\tt{eig}(\cdot,\cdot)$ denotes the eigenvalues of the matrix pencil $\lambda \bE_r^{(k)} - \bE_r^{(k)}$. The algorithm continues until convergence. Upon convergence, the reduced model is guaranteed to be asymptotically stable and satisfies a restricted set of $\Htwo$-optimality conditions. We refer the reader to the original source~\cite{gugercin2008iterative} for details. 
A brief sketch of \SK is given below. 
\begin{algorithm}[htp] 
 \caption{The iterative SVD-rational Krylov (\SK) \cite{gugercin2008iterative}}  
\label{ALG1}                                     
\algorithmicrequire~LTI system described through a state-space realization $(\bA,\bE,\bB,\bC)$, reduced order $r$, an initial selection of shift values $\{\eta_1^{(0)},\eta_2^{(0)},\ldots \eta_r^{(0)}\}$, and a tolerance value $\tau >0$.

Compute projection matrices:
\begin{align}
    \bV_r^{(0)} &= \begin{bmatrix}
        (\eta_1^{(0)} \bE -\bA)^{-1} \bB &  \cdots &  (\eta_r^{(0)} \bE -\bA)^{-1} \bB \nonumber
        \end{bmatrix}, \\
        \bW_r^{(0)} &=  \bQ \bE \bV_r^{(0)}.
\end{align}
Compute reduced-order matrices at step $k=0$:
\begin{align}
    \bA_r^{(0)} &= \left[\bW_r^{(0)}\right]^T \bA \bV_r^{(0)}, \ \ \bE_r^{(0)} = \left[\bW_r^{(0)}\right]^T \bE \bV_r^{(0)},\\
        \bB_r^{(0)} &= \left[\bW_r^{(0)}\right]^T \bB, \ \ \bC_r^{(0)} =  \bC \bV_r^{(0)}.
\end{align}

\algorithmicensure~A reduced-order system given by:  $ \bAtil_r \in \mathbb{R}^{r \times r}, \ \bEtil_r \in \mathbb{R}^{r \times r}, \  \bBtil_r,  \bCtil_r^T \in \mathbb{R}^r$.

\begin{algorithmic} [1]                                        
\WHILE {(the relative change in $\eta_i^{(k)}$) $\geq \tau$}
\STATE Compute matrices $ \bA_r^{(k)}$ and $\bE_r^{(k)}$ as in (\ref{red_mat_stepk_alg1}).
\STATE Assign the following new shifts:
\begin{align}
\eta_i^{(k+1)} \leftarrow -\tt{eig}(\bA_r^{(k)},   \bE_r^{(k)})
\end{align}
\vspace{-4mm}
\STATE $k = k+1$;
\STATE Compute projection matrices as in (\ref{Vrk_alg1}) and in (\ref{Wrk_alg1}).
\STATE Compute matrices at step $k\geq 1$ as in (\ref{red_mat_stepk_alg1}).
\ENDWHILE
\end{algorithmic}
\end{algorithm}
Note that in Algorithm \ref{ALG1}, we have included the $\bE$ matrix explicitly, whereas, in the original contribution \cite{gugercin2008iterative}, $\bE$ was assumed to be identity.


\section{The proposed data-driven reformulation}
\label{sec:main}

Inspired by \cite{gosea2022data} where a fully data-driven formulation of \BT was introduced, the main idea is to produce a data-driven formulation of \SK  using an approximate formulation of the square factor $\bL$ of $\bQ$ in \eqref{lyapFact}, using numerical quadrature. But, the quadrature approximations of $\bL$ and $\bQ$ are \emph{never computed explicitly} since this will indeed require a state-space representation. What is actually needed is the particular structure of matrices $\quadL$, i.e., the approximate square factors, and how it incorporates the input-output data. In the end, the whole formulation will be realization-free and will only necessitate transfer function samples. In the current work, we implicitly apply a quadrature approximation to the Gramian $\bQ$, i.e., $\quadQ \approx \bQ$ (provided that the quadrature scheme is rich enough) where
\begin{align} \label{quad_Q}
\begin{split}
\hspace{-2mm} \quadQ = \sum_{k=1}^{\nq} \phi_k^2  (-\imunit \omega_k \bE^T -\bA^T)^{-1} \bC^T \bC (\imunit \omega_k \bE -\bA)^{-1}.
\end{split}
\end{align}
Here, $\phi_k^2$ and $\imunit \omega_k$ are, respectively, the quadrature weights and nodes. As noted, we will not use $\quadQ$ explicitly but exploit the structure of the approximate square factor
$\quadL$:
\begin{equation}  \label{quad_L}
\quadL^* = \left[ \begin{matrix}
\phi_1 \bC (\imunit \omega_1\bE -\bA)^{-1} \\
\vdots \\ \phi_\nq \bC (\imunit \omega_\nq \bE -\bA)^{-1} 
\end{matrix}  \right] \in \mathbb{C}^{\nq \times n},~~
\quadQ = \quadL \quadL^*.
\end{equation}
As will be shown next, the matrix $\quadL^*$ is not required either, but only as products with other matrices. These matrix products will be interpreted as \emph{data}. The following result directly follows from~\eqref{eq:K} and~\eqref{quad_L}. 
\begin{corollary}
Let $\cK$ and $\quadL$ as defined in \eqref{eq:K} and~\eqref{quad_L}, respectively. Then,
\begin{equation}
    \quadL^* =  \bPhi \cK^T(\bA^T,\bE^T,\bC^T,\imunit \bomega),
\end{equation}
where  $\bPhi = \text{diag}(\phi_1,\phi_2,\ldots,\phi_\nq)$ is a diagonal matrix of quadrature weights.
\end{corollary}

Now, for a particular choice of shifts $\eta_1^{(k)},\eta_2^{(k)},\ldots \eta_r^{(k)}$ associated to step $k$ of the \SK iteration, we can explicitly write the matrix $\bV_r^{(k)} \in \IC^{n \times r}$,  as in Section \ref{sec:ISRK}:
\begin{align}
\begin{split}
    &\bV_r^{(k)}
    = \cK(\bA,\bE,\bB,\eta^{(k)})  \\
    &= \begin{bmatrix}
        (\eta_1^{(k)} \bE -\bA)^{-1} \bB &  \cdots &  (\eta_r^{(k)} \bE -\bA)^{-1} \bB.
    \end{bmatrix}
    \end{split}
\end{align}
Then $\quadW_r^{(k)} \in \IC^{n \times r}$ is constructed  using $\bV_r^{(k)}$ and the quadrature approximation of $\bQ$, i.e. $\quadQ$, as
\begin{equation}\label{eq:quadWrk}
    \quadW_r^{(k)} = \quadQ \bE \bV_r^{(k)} =  \quadL  \quadL^* \bE \bV_r^{(k)}.
\end{equation}
Basically, for the quadrature-based \SK method proposed here, the matrix $\bW_r^{(k)}$ in \eqref{red_mat_stepk_alg1} will be replaced by  $\quadW_r^{(k)}$. We will show that the matrices that compose the ROM at every iteration step, may be written explicitly in terms of transfer function evaluations.

\subsection{Computation of the reduced-order matrices}

The ROM matrices at  iteration step $k$ of \SK (associated with the shifts $\eta_1^{(k)},\eta_2^{(k)},\ldots \eta_r^{(k)}$)  are computed by means of a double-sided projection as 
\begin{align}\label{mat_stepk_new}
\begin{split}
        \bAtil_r^{(k)} &= \left[\quadW_r^{(k)} \right]^* \bA \bV_r^{(k)}, \ \
           \bEtil_r^{(k)} = \left[\quadW_r^{(k)} \right]^* \bE \bV_r^{(k)}, \\
                \bBtil_r^{(k)} &= \left[\quadW_r^{(k)} \right]^* \bB,  \ \  \bCtil_r^{(k)} = \bC \bV_r^{(k)}.
\end{split}
\end{align}

Next, we introduce special matrices that will enter the formulation of the ROM's matrices in \eqref{mat_stepk_new}.
\begin{definition}\label{def:LMVW}
Let matrices $\quadLL^{(k)}, \quadMM^{(k)} \in \IC^{\nq \times r}$ be defined as follows
    \begin{equation}\label{def_2data_mat}
  \quadLL^{(k)} = \quadL^* \bE  \bV_r^{(k)}, \ \  \quadMM^{(k)} = \quadL^*  \bA \bV_r^{(k)}.
\end{equation}
Additionally, let the following vectors be given by
\begin{equation}\label{VW_def}
    \quadVV = \quadL^* \bB \in \IC^{\nq \times 1}, \ \ \quadWW = \bC \bV_r^{(k)} \in \IC^{1 \times r}.
\end{equation}
\end{definition}

Similar to the results in \cite{gosea2022data},  we show that matrices $  \quadLL^{(k)}$ and $ \quadMM^{(k)}$ can be put together in terms of data (samples of the transfer function in \eqref{TF_def}), without requiring access to the original system realization $(\bE,\bA,\bB,\bC)$. The same holds for vectors in \eqref{VW_def}.

\begin{lemma}\label{lem:data_mat}
    The following results hold for any values of $1 \leq i \leq \nq$ and $1 \leq j \leq r$:
    \begin{align}
       \left(\quadLL^{(k)}\right)_{i,j} &=  -\phi_i \frac{\bH(\imunit \omega_i)-\bH(\eta_j^{(k)})}{\imunit \omega_i-\eta_j^{(k)}}, \\
         \left(\quadMM^{(k)}\right)_{i,j} &= -\phi_i \frac{\imunit \omega_i \bH(\imunit \omega_i)-\eta_j^{(k)}\bH(\eta_j^{(k)})}{\imunit \omega_i-\eta_j^{(k)}}, \\
             \left( \quadVV \right)_{i} &= \phi_i \bH(\imunit \omega_i), \ \
     \left(  \quadWW \right)_j = \bH(\eta^{(k)}_j ).
    \end{align}
\end{lemma}

\emph{Proof:} See Section \ref{app:proof_lem} (Appendix).   \unskip\nobreak\hfill $\square$

\noindent
This shows that the matrices defined in (\ref{def_2data_mat}) are scaled Loewner matrices (as used in \cite{mayo2007fsg}), multiplied to the left with the diagonal matrix of quadrature nodes. 

\begin{theorem}The matrices in \eqref{mat_stepk_new} corresponding to the ROM computed at the $k$th step of the proposed iterative procedure are expressed solely in terms of matrices in (\ref{def_2data_mat}), and of vectors in (\ref{VW_def}). More explicitly, we have
\begin{align}\label{data_mat_expr}
\begin{split}
       \bAtil_r^{(k)} &= \left[\quadLL^{(k)}\right]^*   \quadMM^{(k)}, \ \ \bEtil_r^{(k)} = \left[\quadLL^{(k)}\right]^*   \quadLL^{(k)}, \\
         \bBtil_r^{(k)} &= \left[\quadLL^{(k)}\right]^* \quadVV, \ \    \bCtil_r^{(k)} = \quadWW.
         \end{split}
\end{align}
\end{theorem}
\emph{Proof:} \ First, we show that the matrix $\bAtil_r^{(k)}$ in \eqref{mat_stepk_new} can be written as product of data matrices from (\ref{def_2data_mat}), by substituting the formulas of $\quadW_r^{(k)}$ in \eqref{eq:quadWrk}, as:
\begin{align}
\begin{split}
    \bAtil_r^{(k)} &= \left[\quadW_r^{(k)} \right]^* \bA \bV_r^{(k)} = \left[ \quadL  \quadL^* \bE \bV_r^{(k)} \right]^* \bA \bV_r^{(k)}\\
    &=  \left[\quadL^* \bE  \bV_r^{(k)}  \right]^* \quadL^*  \bA \bV_r^{(k)} =   \left[\quadLL^{(k)}\right]^*   \quadMM^{(k)},
    \end{split}
\end{align}
and similarly for the matrix $\bE_r^{(k)}$:
\begin{align}
\begin{split}
    \bEtil_r^{(k)} &= \left[\quadW_r^{(k)} \right]^* \bE \bV_r^{(k)} = \left[ \quadL  \quadL^* \bE \bV_r^{(k)} \right]^* \bE \bV_r^{(k)}\\
    &=  \left[\quadL^* \bE  \bV_r^{(k)}  \right]^* \quadL^*  \bE \bV_r^{(k)} =   \left[\quadLL^{(k)}\right]^*   \quadLL^{(k)}.
    \end{split}
\end{align}
Then, we write the row and column vectors, $\bBtil_r^{(k)}$ and $\bCtil_r^{(k)}$, in terms of data vectors introduced in (\ref{VW_def}), as follows:
\begin{align}
\begin{split}
    \bBtil_r^{(k)} &= \left[\quadW_r^{(k)} \right]^* \bB = \left[ \quadL  \quadL^* \bE \bV_r^{(k)} \right]^* \bB \\
    &= \left[\quadL^* \bE  \bV_r^{(k)}  \right]^* \quadL^*  \bB = \left[\quadLL^{(k)}\right]^* \quadVV, \\
     \bCtil_r^{(k)} &= \bC \bV_r^{(k)} = \quadWW.  \qed
    \end{split}
    \vspace{-5mm}
\end{align}

This result is of importance since it illustrates how to explicitly compute the ROM matrices~\eqref{mat_stepk_new} in \SK, simply utilizing data-driven components only (put together by using transfer function evaluations, quadrature nodes, shifts, and weights). This data-driven construction structure is similar to~\cite{gosea2022data}. However, there are some key differences: Unlike \QBT~\cite{gosea2022data}, for the proposed method here, one performs an iteration (to compute ``sub-optimal" shifts upon convergence) and decides on a truncation order $r$, a priori. Actually, the method under study here resembles TF-\RK in \cite{BeattieGugercin2012RlznIndH2Approx}. The latter is also a purely data-driven approach and relies on an iteration based on Petrov-Galerkin rational Krylov projection computed via data-based Loewner matrices.

\begin{remark}
We note that if $r = \nq$, and that if at step $k$ the matrix $\quadLL^{(k)}$ is invertible, then the realization in \eqref{data_mat_expr}, can be equivalently rewritten as:
\begin{align}
\begin{split}
    \hspace{-2mm}   \bAtil_r^{(k)} &=   \quadMM^{(k)}, \ \ \bEtil_r^{(k)} =  \quadLL^{(k)}, \ \
         \bBtil_r^{(k)} = \quadVV, \ \    \bCtil_r^{(k)} = \quadWW,
         \end{split}
\end{align}
which is precisely the way one would write a realization in the Loewner framework \cite{mayo2007fsg}. This is obtained by multiplying (to the left) all system matrices appearing in the differential state equation from \eqref{data_mat_expr} with the inverse of matrix $\left[\quadLL^{(k)}\right]^* \in \IC^{r \times r}$. 

In practice, one would expect that the reduction order $r$ is much smaller than the number of quadrature nodes $\nq$. A similar result can be derived in this more realistic case, i.e.,   $r \ll \nq$, as shown above. One would then need to use instead the right generalized inverse of matrix $\left[\quadLL^{(k)}\right]^* \in \IC^{r \times \nq}$.
\end{remark}

\subsection{The proposed data-driven formulation of \SK}

Now based on the theory developed earlier, we present the proposed data-driven algorithm. It is an iterative procedure, which terminates whenever the difference between shifts in neighboring steps is small enough (based on a pre-defined tolerance value $\tau >0$).  The brief sketch of the proposed method, Quadrature-based (data-driven) iterative SVD-rational Krylov Algorithm (\QSK), is given in Algorithm~\ref{ALG2} below. As it is clear from the algorithm, at no point, one would need a state-space representation. \QSK only uses transfer function samples. 

\begin{algorithm}[htp] 
 \caption{Quadrature-based (data-driven) iterative SVD-rational Krylov Algorithm (\QSK)}  
\label{ALG2}                                     
\algorithmicrequire~LTI system described through a transfer function evaluation map, $\bH(s)$, quadrature nodes, $\omega_k$, and weights, $\phi_k$, for $k=1,2,\ldots,\nq$, a truncation value given by $1\leq r\leq \nq$, and an initial selection of shift values $\{\eta_1^{(0)},\eta_2^{(0)},\ldots \eta_r^{(0)}\}$, and a tolerance value $\tau >0$. 

\algorithmicensure~A reduced-order system given by:  $ \bAtil_r \in \mathbb{R}^{r \times r}, \ \bEtil_r \in \mathbb{R}^{r \times r}, \  \bBtil_r,  \bCtil_r^T \in \mathbb{R}^r$.

\begin{algorithmic} [1]                                        
\WHILE {(the relative change in $\eta_i^{(k)}$) $\geq \tau$}
\STATE \label{sample} Sample transfer function values $\{H(\imunit \omega_i)\}_{i=1}^{\nq}$ and $\{H(\eta_j^{(k)})\}_{j=1}^{r}$.
\STATE Based on the data values (samples of the transfer function) and quadrature weights $\{\phi_k\}$, put together the data matrices in (\ref{def_2data_mat}) and data vectors in (\ref{VW_def}).
	\STATE Construct the matrices of the reduced-order system as in (\ref{data_mat_expr}), i.e., $\bAtil_r^{(k)}, \bEtil_r^{(k)}, \bBtil_r^{(k)}$, and $\bCtil_r^{(k)}$.
	\STATE Assign the following new shifts:
\begin{align}
\eta_i^{(k+1)} \leftarrow -\tt{eig}(\bA_r^{(k)},   \bE_r^{(k)})
\end{align}
\ENDWHILE
\STATE $\bAtil_r \leftarrow \bA_r^{(k)}$, $\bEtil_r \leftarrow \bE_r^{(k)}$, $\bBtil_r \leftarrow \bB_r^{(k)}$, $\bCtil_r \leftarrow \bC_r^{(k)}$.
\end{algorithmic}
\end{algorithm}

The original \RK  in \cite{gugercin2008h_2} for $\Htwo$-optimal model reduction typically shows rapid convergence in a wide range of examples; however, proof of convergence, in general, is not guaranteed except for the special cases of symmetric systems, see \cite{flagg2012convergence} where \RK was proven to be locally convergent. \SK has a similar iterative structure and with the inclusion of the observability Gramian $\mathbf{Q}$, it shows even a better convergence behavior, see~\cite{gugercin2008iterative}.
Since the iteration in \QSK is inspired by \SK, we believe that it will have similar convergence properties, which is indeed the case in the numerical example below. However, a detailed analysis could be necessary to understand the impact of the quadrature approximation on its convergence property. 


\section{Numerical example: A classical MOR benchmark}
\label{sec:num}

We illustrate the performance of the realization-free method presented in this work, i.e., \QSK, on a MOR benchmark example: the clamped beam model from \cite{morChaV02}. The original system dimension is $n=348$, with $m=p=1$. We compare \QSK to  \RK in \cite{gugercin2008h_2} and the original projection-based \SK in \cite{gugercin2008iterative}. For \QSK developed in this work, we have used $\nq = 400$ purely imaginary quadrature nodes ($200$ in the interval $[10^{-2},10^2]\imunit$ and their complex conjugates). This is necessary to enforce ROMs with a real-valued realization. The quadrature scheme used here for the implicit approximation of the $\bQ$ Gramian in \eqref{quad_Q}, is the classical exponentially-convergent trapezoidal scheme (see \cite{gosea2022data} for more details). The tolerance value, which acts as a stopping criterion, is chosen as $\tau = 10^{-4}$.

\begin{figure}[!ht]
\includegraphics[width = 0.47\textwidth]{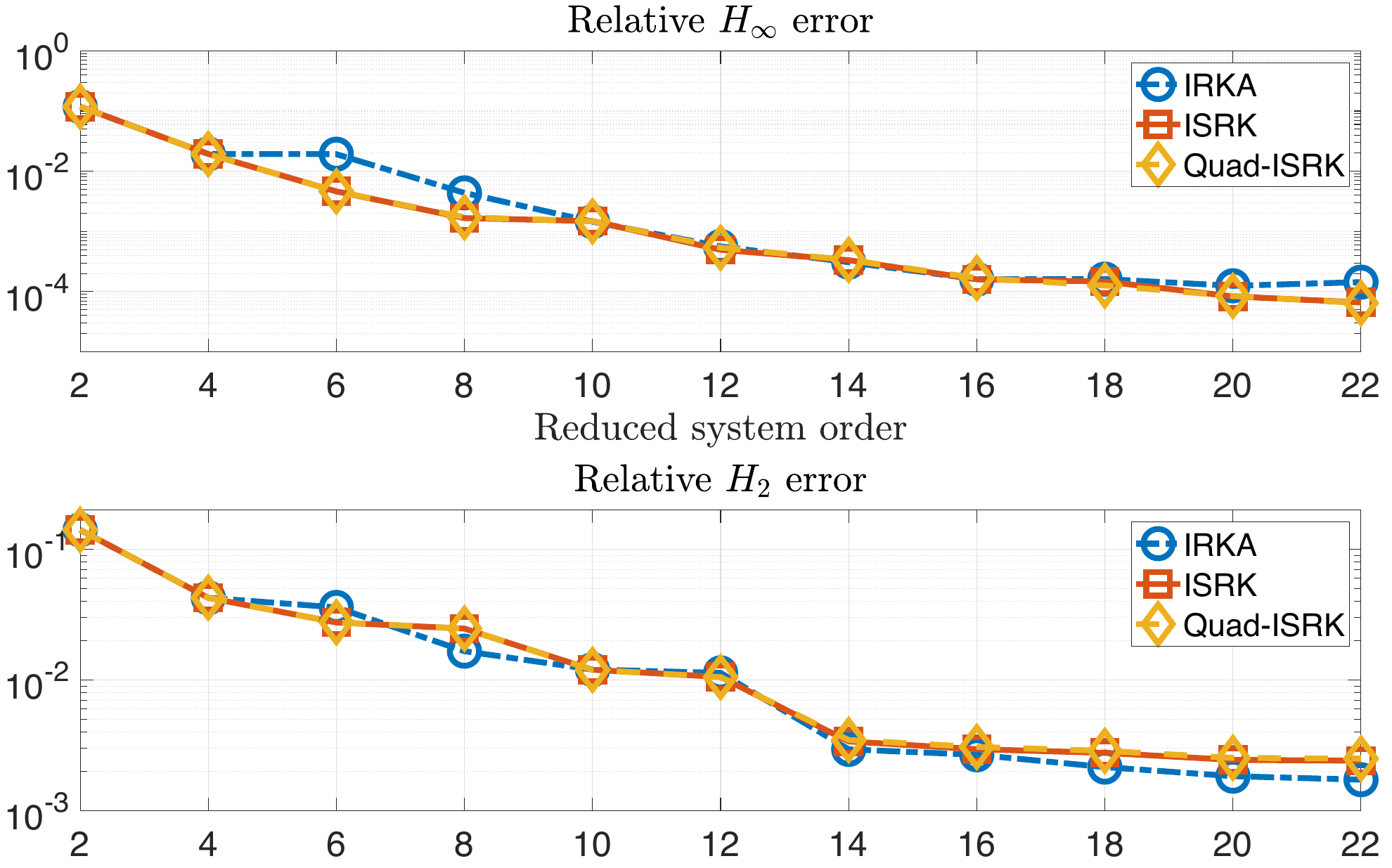}
\vspace{-2mm}
      \caption{Relative $\mathcal{H}_2$ norms of the error systems.}
	\label{fig:iss12}
\end{figure}

The top plot in Figure~\ref{fig:iss12} shows the relative $\Hinf$ error for all three methods as $r$ varies. As the figure shows, \QSK almost identically follows the behavior of the original projection-based formulation \SK; thus we match the quality of \SK by only using transfer function samples. The situation is similar in the bottom plot in Figure~\ref{fig:iss12} where we depict the $\Htwo$-error vs reduced order. Once again, the data-driven \QSK is able to replicate the behavior of \SK. As expected, for many $r$ values, \RK has a smaller $\Htwo$ error compared to \QSK and \SK as it is explicitly tailored towards the $\Htwo$ norm.

\section{Conclusions and future work}
\label{sec:conc}

We have developed~\QSK, a data-driven formulation of \SK, originally an intrusive projection-based method. As in quadrature-based balancing approaches, the idea is to approximate the original system's infinite Gramian with a finite quadrature scheme. The numerical quadrature formulations are never explicitly computed. We showed that at each step of the proposed iteration, one does not require access to the original system's matrices, but only to transfer function samples evaluated at particular values given by quadrature nodes (on the imaginary axis), or the rational Krylov shifts computed in the previous step. The numerical example shows that the proposed data-driven formulation~\QSK almost exactly replicates the behavior of the original projection-based one~\SK, thus illustrating its efficiency.  Future endeavors may include investigating the effect of choosing different quadratures, such as adaptive schemes, and also extending the new method to other classes of systems (structured linear or even mildly nonlinear).

\section{Appendix: proof of Lemma \ref{lem:data_mat}}
\label{app:proof_lem}

To show that the matrices introduced in Definition \ref{def:LMVW} are indeed purely data-based, we start by noting that the $(i,j)$th element of matrix  $\quadLL^{(k)} \in \IC^{\nq \times r}$ can be written as:
\begin{align}\label{ij_Loew}
\begin{split}
   \bfe_i^T \quadLL^{(k)} \bfe_j &=  \bfe_i^T \quadL^* \bE  \bV_r^{(k)} \bfe_j\\
   &= \phi_i \bC (\imunit \omega_i\bE -\bA)^{-1} \bE (\eta_j^{(k)} \bE -\bA)^{-1} \bB,
    \end{split}
\end{align}
where $\bfe_\ell = \begin{bmatrix} 0 & 0 & \ldots & 0 & 1 & 0 & \ldots & 0  \end{bmatrix} \in \IR^n$ is the $\ell$th unit vector of length $n$. Now, by making use of the simple identity also utilized in the Loewner framework, i.e.,
\begin{equation}
    \bE =  \frac{(x \bE - \bA) -  (y \bE - \bA)}{x-y}, 
\end{equation}
for all $x,y \in \IC$, and $\bA, \bE \in \IC^n$, it follows that:
\begin{equation*}
    (x \bE - \bA)^{-1} \bE  (y \bE - \bA)^{-1}  =  -\frac{(x \bE - \bA)^{-1} - (y \bE - \bA)^{-1}}{x-y}.
\end{equation*}
Next, by substituting  $x = \imunit \omega_i$ and $y = \eta_j^{(k)}$ in the identity above and then into (\ref{ij_Loew}), it follows that:
\begin{align}
\begin{split}
   \bfe_i^T \quadLL \bfe_j &= -\phi_i \bC \frac{(\imunit \omega_i \bE - \bA)^{-1} - (\eta_j^{(k)} \bE - \bA)^{-1}}{\imunit \omega_i-\eta_j^{(k)}} \bB,
     \end{split}
\end{align}
which proves the first result.
Similarly, we observe  that the $(i,j)$th element of matrix  $\quadMM^{(k)} \in \IC^{\nq \times r}$ can be written as:
\begin{align}\label{ij_sh_Loew}
\begin{split}
   \bfe_i^T \quadMM^{(k)} \bfe_j &=  \bfe_i^T \quadL^* \bA  \bV_r^{(k)} \bfe_j\\
   &= \phi_i \bC (\imunit \omega_i\bE -\bA)^{-1} \bA (\eta_j^{(k)} \bE -\bA)^{-1} \bB.
    \end{split}
\end{align}
Now, by making use of another standard equality, i.e.,
\begin{equation}
    \bA =  \frac{y(x \bE - \bA) -  x(y \bE - \bA)}{x-y}, 
\end{equation}
for all $x,y \in \IC$, and $\bA, \bE \in \IC^n$, it follows that:
\begin{equation*}
    (x \bE - \bA)^{-1} \bA  (y \bE - \bA)^{-1}  =  -\frac{x(x \bE - \bA)^{-1} - y(y \bE - \bA)^{-1}}{x-y}.
\end{equation*}
Next, by substituting in the identity above $x = \imunit \omega_i$ and $y = \eta_j^{(k)}$, and then into (\ref{ij_Loew}), it follows that:
\begin{align}
   \bfe_i^T \quadMM \bfe_j &= -\phi_i \bC \frac{\imunit \omega_i (\imunit \omega_i \bE - \bA)^{-1} - \eta_j^{(k)} (\eta_j^{(k)} \bE - \bA)^{-1}}{\imunit \omega_i-\eta_j^{(k)}} \bB \nonumber
\end{align}
which proves the second result.
Finally, the $i$th entry of $\quadVV$, and  the $j$th entry of $\quadWW$ are written as:
\begin{align*}
    \bfe_i^T \quadVV &= \phi_i \bC (\imunit \omega_i\bE -\bA)^{-1} \bB = \phi_i \bH(\imunit \omega_i), \\
    \quadWW \bfe_j &= \bC (\eta^{(k)}_j \bE -\bA)^{-1} \bB = \bH(\eta^{(k)}_j ).
\end{align*}

\bibliography{literature}    

\end{document}